\documentclass[reqno,dvips,12pt]{article} %% english should be theLAST
\usepackage{epsfig,graphicx,amsmath,amssymb,amsfonts}
cm \textwidth=16.0 true cm \emergencystretch=10pt

\newtheorem{theorem}{\qquad Theorem}
\newtheorem{lemma}{\qquad Lemma} %[section]
\newcommand{\abs}[1]{\left\vert{#1}\right\vert}
\newcommand{\e}{\varepsilon}

\DeclareMathOperator{\re}{re}

\newcommand{\cA}{{\cal A}}
\newcommand{\cM}{{\cal M}}

\newcommand{\cN}{{\cal N}}
\newcommand{\cP}{{\cal P}}

\newcommand{\nn}{\nonumber}

\newcommand{\pa}{\partial}

\begin{document}

\title {\textbf{On a Distribution in Frequency Probability Theory
Corresponding to the Bose-Einstein Distribution}}

\author{\textbf{V.P.Maslov}\thanks{Moscow Institute of Electronics and Mathematics,
pm@miem.edu.ru}}
\date{ }
%\emph{Moscow Institute of Electronics and Mathematics, pm@miem.edu.ru}

\maketitle
\begin{abstract}
The notion of density of a finite set is discussed. We proof a
general theorem of set theory which refines Bose-Einstein
distribution.
\end{abstract}

In our paper \cite{General_theor} a mistake through negligence
occurred\footnote{On page 7 at reference to the paper
\cite{NelinSred_ArXiv} instead of formula $(95)$ we should refer
to formula $(94)$ and consider $(96)$, namely $
\nu=\nu'+O(\frac{1}{N^{1/4}})$. This changes the evaluation $(4)$:
the exponent should be $3/4$ instead of $1/2$.}. Here we produce
a detailed proof of the theorem under consideration.

As in \cite{{NelinSred_ArXiv}}, the values of the random variable
$ x_1, \dots,  x_s$ are ordered in absolute value. Some of the
numbers $ x_1,\dots, x_s$ may coincide. Then these numbers are
combined adding the corresponding ''probabilities'', i.e., the
ratio of the number of ``hits'' at~$ x_i$ to the general number
of trials. The number of equal $ x_i:  x_i= x_{i+1}= \dots =
x_{i+k}$ will be called the multiplicity $q_i$ of the value $
x_i$. In our consideration, both the number of trials $N$ and~$s$
tend to infinity \cite{Teorver}.

Let $N_i$ be the number of ''appearances'' of the value $ x_i: \
 x_i <  x_{i+1}$, then
\begin{equation}
\sum^s_{i=1} \frac{N_i}{N}  x_i=M, \label{Zipf1}
\end{equation}
where $M$ is the mathematical expectation.

The cumulative probability $\cP_k$ is the sum of the first~$k$
probabilities in the sequence $ x_i$: $\cP_k=\frac 1N
\sum_{i=1}^k N_i$, where $k<s$. We denote $NP_k=B_k$.

If all the variants for which
\begin{equation}\label{A}
\sum_{i=1}^s N_i = N
\end{equation}
and
\begin{equation}\label{B}
\sum_{i=1}^s N_i  x_i \leq E, \ \ E=MN\leq N \overline{ x},
\end{equation}
where $\overline{ x}=\frac{\sum_{i=1}^s q_i  x_i}{Q}$,
$Q=\sum_{i=1}^s q_i$, are equivalent (equiprobable), then
\cite{MatZamTheor,NelinSred} the majority of the variants will
accumulate near the following dependence of the ''cumulative
probability'' $B_l\{N_i\}=\sum_{i=1}^l N_i$,
\begin{equation}
\sum_{i=1}^l N_i= \sum_{i=1}^l \frac{q_i}{e^{\beta' x_i-\nu'}-1},
\label{Zipf2}
\end{equation}
where $\beta'$ and $\nu'$ are determined by the conditions
\begin{equation}\label{Zipf2a}
B_s=N,
\end{equation}
\begin{equation}\label{Zipf2a'}
\sum_{i=1}^s \frac{q_i  x_i}{e^{\beta'  x_i-\nu'}-1}=E,
\end{equation}
as $N \to \infty$ and $s \to \infty$.

We introduce the notation: $\cM$ is the set of all sets $\{N_i\}$
satisfying conditions~(\ref{A}) and~(\ref{B}); \ $\cN\{\cM\}$ is
the number of elements of the set~$\cM$.

\begin{theorem} \label{theor1}
Suppose that all the variants of sets $\{N_i\}$ satisfying the
conditions ~(\ref{A}) and ~(\ref{B}) are equiprobable. Then the
number of variants $\cN$ of sets $\{N_i\}$ satisfying
conditions~(\ref{A}) and~(\ref{B}) and the additional relation
\begin{equation} \label{theorema1}
|\sum^l_{i=1} N_i - \sum^l_1\frac{q_i}{e^{\beta'
 x_i-\nu'}-1}|\geq  N^{(3/4+\varepsilon)}
\end{equation}
is less than $\frac{c_1 \cN\{\cM\}}{N^m}$
(where~$c_1$ and~$m$ are any arbitrary numbers,
$\sum_{i=1}^l q_i \geq\varepsilon Q$, and
$\varepsilon$ is arbitrarily small).
\end{theorem}

{\bf{\qquad Proof of Theorem 1.}}

Let $\cA$ be a subset of $\cM$ satisfying the condition
$$
|\sum_{i=l+1}^s N_i - \sum_{i=l+1}^s \frac{q_i} {e^{\beta
x_i-\nu}-1}|\leq \Delta;
$$
$$
|\sum_{i=1}^l N_i-\sum_{i=1}^l \frac{q_i} {e^{\beta'
x_i-\nu'}-1}|\leq \Delta,
$$
where $\Delta$, $\beta$, $\nu$ are some real numbers independent of~$l$.

We denote
$$
|\sum_{i=l+1}^s N_i-\sum_{i=l+1}^s \frac{q_i} {e^{\beta
x_i-\nu}-1}| =S_{s-l};
$$
$$
|\sum_{i=1}^l N_i-\sum_{i=1}^l \frac{q_i} {e^{\beta'
x_i-\nu'}-1}| =S_l.
$$

Obviously, if $\{N_i\}$ is the set of all sets of integers on the
whole, then
\begin{equation}\label{Proof1}
\cN\{\cM \setminus \cA\} = \sum_{\{N_i\}} \Bigl(
\Theta(E-\sum_{i=1}^s N_i x_i) \delta_{(\sum_{i=1}^s N_i),N}
\Theta(S_l-\Delta)\Theta (S_{s-l}-\Delta)\Bigr),
\end{equation}
where $\sum N_i=N$.

Here the sum is taken over all integers $N_i$, $\Theta( x)$ is
the Heaviside function, and $\delta_{k_1,k_2}$ is the Kronecker
symbol.

We use the integral representations
\begin{eqnarray}
&&\delta_{NN'}=\frac{e^{-\nu N}}{2\pi}\int_{-\pi}^\pi d\varphi
e^{-iN\varphi} e^{\nu N'}e^{i N'\varphi},\label{D7}\\
&&\Theta(y)=\frac1{2\pi i}\int_{-\infty}^\infty d x\frac1{
x-i}e^{\beta y(1+i x)}.\label{D8}
\end{eqnarray}

Now we perform the standard regularization.
We replace the first Heaviside function~$\Theta$
in~(\ref{Proof1}) by the continuous function
\begin{equation}
\Theta_{\alpha}(y) =\left\{
\begin{array}{ccc}
0 &\mbox{for}& \alpha > 1, \  y<0 \nn \\
1-e^{\beta y(1-\alpha)} &\mbox{for}& \alpha > 1,\  y \geq 0,
\label{Naz1}
\end{array}\right.
\end{equation}
\begin{equation}
\Theta_{\alpha}(y) =\left\{
\begin{array}{ccc}
e^{\beta y(1-\alpha)} &\mbox{for}&\alpha < 0, \ y<0 \nn \\
1 &\mbox{for}& \alpha < 0, \ y \geq 0, \label{Naz2}
\end{array}\right.
\end{equation}
where $\alpha \in (-\infty,0) \cup (1, \infty)$ is a parameter,
and obtain
\begin{equation}\label{proof2}
\Theta_\alpha(y) = \frac1{2\pi i} \int_{-\infty}^{\infty}
e^{\beta y(1+ix)} (\frac 1{x-i} - \frac 1{x-\alpha i}) dx.
\end{equation}

If  $\alpha > 1$, then $\Theta(y)\leq \Theta_{\alpha}(y)$.

Let $\nu <0$. We substitute~(\ref{D7}) and~(\ref{D8})
into~(\ref{Proof1}),  interchange the integration and summation,
then pass to the limit as $\alpha \to \infty$
and obtain the estimate
\begin{eqnarray}
&&\cN\{\cM \setminus   \cA\} \leq \nn \\
&&\leq \Bigl|\frac{e^{-\nu N+\beta E }}{i(2\pi)^2}\int_{-\pi}^\pi
\bigl[ \exp(-iN\varphi)
\sum_{\{N_j\}}\Bigl(\exp\bigl\{(-\beta\sum_{j=1}^s N_j
x_j)+(i\varphi+\nu)
N_j\bigr\}\bigr]\ d\varphi \times \nn \\
&& \times \Theta(S_l -\Delta)\Theta(S_{s-l}-\Delta)\Bigr|, \ \sum
N_i=N,
\end{eqnarray}
where $\beta$ and $\nu$ are real parameters
such that the series converges for them.

To estimate the expression in the right-hand side,
we bring the absolute value sign inside the integral sign
and then inside the sum sign, integrate over $\varphi$, and obtain
\begin{eqnarray}
&&\cN\{\cM \setminus \cA\} \leq \frac{e^{-\nu N+\beta E }}{2\pi}
\sum_{\{N_i\}}\exp\{-\beta\sum_{i=1}^sN_i x_i+\nu
N_i\}\times \nn \\
&& \times\Theta (S_l-\Delta)\Theta (S_{s-l}-\Delta).
\end{eqnarray}

We denote
\begin{equation}\label{D9a}
Z(\beta,N)=\sum_{\{N_i\}} e^{-\beta\sum_{i=1}^s N_i x_i},
\end{equation}
where the sum is taken over all $N_i$ such that $\sum_{i=1}^s N_i=N$,
$$
\zeta_l(\nu,\beta)= \prod_{i=1}^{l} \xi_i\left(\nu,\beta\right);
\zeta_{s-l}(\nu,\beta)= \prod_{i=l+1}^{s}
\xi_i\left(\nu,\beta\right);
$$
$$
\quad \xi_i(\nu,\beta)= \frac{1}{(1-e^{\nu-\beta x_i})^{q_i}},
\qquad i=1,\dots,l.
$$

It follows from the inequality for the hyperbolic cosine
$\cosh(x)=(e^x+e^{-x})/2$ for $|x_1| \geq \delta; |x_2| \geq
\delta$:
\begin{equation}
\cosh(x_1)\cosh(x_2)= \cosh(x_1+x_2) + \cosh(x_1 - x_2) >
\frac{e^\delta}{2} \label{D33}
\end{equation}
that the inequality
\begin{equation}
\Theta(S_{s-l}-\Delta) \Theta(S_{l}-\Delta)\le e^{-c\Delta}
\cosh\Bigl(c\sum_{i=1}^{l} N_i
-c\phi_l\Bigr)\cosh\Bigl(c\sum_{i=l+1}^{s} N_i
-c\overline{\phi}_{s-l}\Bigr), \label{D34}
\end{equation}
where
$$
\phi_l= \sum_{i=1}^l \frac{q_i}{e^{\beta' x_i-\nu'}-1}; \qquad
\overline{\phi}_{s-l}= \sum_{i=l+1}^s \frac{q_i}{e^{\beta
x_i-\nu}-1},
$$
holds for all positive $c$ and~$\Delta$.

We obtain
\begin{eqnarray}
&&\cN\{\cM \setminus \cA\} \leq  e^{-c\Delta} \exp\left(\beta
E-\nu N\right) \times \nn \\
&& \times \sum_{\{N_i\}}\exp\{-\beta\sum_{i=1}^l N_i
x_i+\nu\sum_{i=1}^l N_i\} \cosh\left(\sum_{i=1}^{l} c N_i -
c\phi\right) \times \nn \\
&& \times \exp\{-\beta\sum_{i=l+1}^s N_i x_i +\nu \sum_{i=l+1}^s
N_i\} \cosh\Bigl(\sum_{i=l+1}^s c N_i -c\overline{\phi}\Bigr) = \nn \\
&& =e^{\beta E} e^{-c\Delta} \times \nn \\
&&\times \left( \zeta_l(\nu-c,\beta) \exp(-c\phi_{l})
+\zeta_l(\nu+c,\beta)\exp(c\phi_{l})\right) \times \nn \\
&&\times\left(\zeta_{s-l}(\nu-c,\beta)
\exp(-c\overline{\phi}_{s-l})
+\zeta_{s-l}(\nu+c,\beta)\exp(c\overline{\phi}_{s-l})\right).
\label{5th}
\end{eqnarray}

Now we use the relations
\begin{equation}\label{5tha}
\frac {\pa}{\pa\nu}\ln \zeta_l|_{\beta=\beta',\nu=\nu'}\equiv
\phi_l; \frac {\pa}{\pa\nu}\ln
\zeta_{s-l}|_{\beta=\beta',\nu=\nu'}\equiv \overline{\phi}_{s-l}
\end{equation}
and the expansion $\zeta_l(\nu\pm c,\beta)$ by the Taylor formula.
There exists a $ \gamma <1$ such that
$$
\ln(\zeta_l(\nu\pm c,\beta)) =\ln\zeta_l(\nu,\beta)\pm
c(\ln\zeta_l)'_\nu(\nu,\beta)+\frac{c^2}{2}(\ln\zeta_l)^{''}_\nu
(\nu\pm\gamma c,\beta).
$$
We substitute this expansion, use formula~(\ref{5tha}), and see that
$\phi_{\nu,\beta}$ is cancelled.

Another representation of the Taylor formula implies
\begin{eqnarray}
&&\ln\left(\zeta_l(\nu+c,\beta)\right)=
\ln\left(\zeta_l(\beta,\nu)\right)+
\frac{c}\beta\frac{\pa}{\pa\nu}\ln\left(\zeta_l(\beta,\nu)\right)+\nn\\
&&+\int_{\nu}^{\nu+c/\beta}d\nu' (\nu+c/\beta-\nu')
\frac{\pa^2}{\pa\nu'^2}\ln\left(\zeta_l(\beta,\nu')\right).\label{CC1}
\end{eqnarray}
A similar expression holds for $\zeta_{s-l}$.

From the explicit form of the function $\zeta_l(\beta,\nu)$, we obtain
\begin{equation}
\frac{\pa^2}{\pa\nu^2}\ln\left(\zeta_l(\beta,\nu)\right)=
\beta^2\sum_{i=1}^{l} \frac{g_i\exp(-\beta(
x_i+\nu))}{(\exp(-\beta( x_i+\nu))-1)^2}\leq \beta^2Qd,
\label{CC2}
\end{equation}
where $d$ is given by the formula
$$
d=\frac{\exp(-\beta( x_1+\nu))}{(\exp(-\beta( x_1+\nu))-1)^2}..
$$
The same estimate holds for $\zeta_{s-l}$.

Taking into account the fact that $\zeta_l\zeta_{s-l}=\zeta_s$,
we obtain the following estimate for $\beta=\beta'$ and $\nu=\nu'$:
\begin{equation}\label{eval1}
\cN\{\cM \setminus \cA\}
\leq\zeta_s(\beta',\nu')\exp(-c\Delta+\frac{c^2}{2}\beta^2Qd)
\exp(E\beta'-\nu'N).
\end{equation}

Now we express $\zeta_s(\nu',\beta')$ in terms $Z(\beta,N)$.
To do this, we prove the following lemma.

\begin{lemma}%lemma 1
Under the above assumptions, the asymptotics of the integral
\begin{equation}\label{lemma_1}
Z(\beta,N) = \frac{e^{-\nu N}}{2\pi}\int_{-\pi}^\pi
 d\alpha e^{-iN\alpha}\zeta_s(\beta,\nu+i\alpha)
\end{equation}
has the form
\begin{equation}\label{lemma_2}
Z(\beta,N) = C e^{-\nu N} \frac{\zeta_s(\beta,\nu)}{|(\partial^2
\ln\zeta_s(\beta,\nu))/ (\partial^2\nu)|} (1+O(\frac 1N)),
\end{equation}
where $C$ is a constant.
\end{lemma}

The proof of the lemma readily follows from the saddle-point method
and the inequalities
\begin{equation}\label{nerav}
|\xi_i(\nu+i\alpha,\beta)| < \xi_i(\nu,\beta), \qquad
|\zeta_s(\nu+i\alpha,\beta)| <\zeta_s(\nu\,\beta),
\end{equation}
which hold, because $e^{\nu-\beta x_i} <1$ for all $\alpha\neq
2\pi n$, where $n$ is an integer. It follows from these
inequalities that $\alpha=0$ is a saddle point of
integral~(\ref{lemma_1})~\cite{Fedoruk}.

We present another proof \cite{TeorVozm}.

%\begin{small} {
We have
\begin{equation}
 Z(\beta,N) = \frac{e^{-\nu N}}{2\pi}\int_{-\pi}^\pi
 e^{-iN\alpha}\zeta_s(\beta,\nu+i\alpha)\,d\alpha
 =\frac{e^{-\nu N}}{2\pi}\int_{-\pi}^\pi e^{NS(\alpha,N)} d\alpha ,\label{D15}
\end{equation}
where
\begin{equation}\label{qq}
    S(\alpha,N) = -i\alpha+ \ln \zeta_s (\beta, \nu +i\alpha)
    = -i\alpha - \sum_{i=1}^s q_i\ln [1-e^{\nu+i\alpha-\beta x_i}].
\end{equation}
Here $S$ depends on $N$, because $s$, $ x_i$, and $\nu$ also
depend on~$N$; the latter is chosen so that the point $\alpha=0$
be a stationary point of the phase~$S$, i.e., from the condition
\begin{equation}\label{qq1}
 N=\sum_{i=1}^s\frac{q_i}{e^{\beta x_i-\nu}-1}.
\end{equation}
We assume that $a_1N \leq s \leq a_2N$, $a_1,
a_2=\operatorname{const}$, and, in addition, $0\le x_i\le B$ and
$B=\operatorname{const}$, $i=1,\dots,s$.
 If these conditions are satisfied
in some interval $\beta\in[0,\beta_0]$ of the values of the inverse temperature,
then all the derivatives of the phase are bounded,
the stationary point is nondegenerate,
and the real part of the phase outside a neighborhood of zero
is strictly less than its value at zero minus some positive number.
Therefore, calculating the asymptotics of the integral,
we can replace the interval of integration $[-\pi,\pi]$
by the interval $[-\e,\e]$.
In this integral, we perform the change of variable
\begin{equation}\label{qqq}
  z=\sqrt {S(0,N)-S(\alpha,N)}.
\end{equation}
This function is holomorphic in the disk $\abs{\alpha}\le\e$
in the complex $\alpha$-plane and has a holomorphic inverse
for a sufficiently small~$\e$.
As a result, we obtain
\begin{equation}\label{qqq1}
    \int_{-\e}^\e e^{NS(\alpha,N)}
    d\alpha=e^{NS(0,N)}\int_\gamma e^{-Nz^2}f(z)\,dz,
\end{equation}
where the path $\gamma$ in the complex $z$-plane is obtained
from the interval $[-\e,\e]$ by the change~\eqref{qqq}
and
\begin{equation}\label{qqq3}
    f(z)=\left(\frac{\partial\sqrt {S(0,N)-S(\alpha,N)}}
    {\partial\alpha}\right)^{-1}\bigg|_{\alpha=\alpha(z)}.
\end{equation}
For a small~$\e$ the path $\gamma$ lies completely inside
the double sector $\re(z^2)>c(\re z)^2$ for some $c>0$;
hence it can be ``shifted'' to the real axis
so that the integral does not change up to terms
that are exponentially small in~$N$.
Thus, with the above accuracy, we have
\begin{equation}\label{qqq4}
  Z(\beta,N) =  \frac{e^{-\nu N}}{2\pi}\int_{-\e}^\e e^{-Nz^2}f(z)\,dz.
\end{equation}
Since the variable $z$ is now real,
we can assume that the function $f(z)$ is finite
(changing it outside the interval of integration),
extend the integral to the entire axis
(which again gives an exponentially small error),
and then calculate the asymptotic expansion of the integral
expanding the integrand in the Taylor series in~$z$ with a remainder.
This justifies that the saddle-point method can be applied
to the above integral in our case.

%} \end{small}

\begin{lemma}%lemma 2
The quantity
\begin{equation}\label{qqq5}
\frac{1}{\cN(\cM)} \sum_{\{N_i\}} e^{-\beta\sum_{i=1}^s N_i x_i},
\end{equation}
where $\sum N_i =N$ and $ x_iN_i\leq E-N^{1/2+\varepsilon}$,
tends to zero faster than $N^{-k}$ for any $k$, $\varepsilon>0$.
\end{lemma}

We consider the point of minimum in $\beta$
of the right-hand side of ~(\ref{5th})
with $\nu(\beta,N)$ satisfying the condition
$$
\sum \frac{q_i}{e^{\beta x_i-\nu(\beta,N)}-1} =N.
$$
It is easy to see that it satisfies condition~(\ref{Zipf2a}).
Now we assume that the assumption of the lemma is not satisfied.

Then for $\sum N_i=N$,  $\sum  x_i N_i\geq E-N^{1/2+\varepsilon}$,
we have
$$
e^{\beta E}\sum_{\{N_i\}} e^{-\beta\sum_{i=1}^s N_i x_i}\geq
e^{(N^{1/2}+\varepsilon)\beta}.
$$
Obviously, $\beta\ll \frac{1}{\sqrt{N}}$
provides a minimum of~(\ref{5th})
if the assumptions of Lemma~1 are satisfied,
which contradicts the assumption that
the minimum in~$\beta$ of the right-hand side of~(\ref{5th})
is equal to~$\beta'$.

We note that, following the same scheme,
G.~V.~Koval' obtained a direct proof of this lemma.

We set $c=\frac\Delta{N^{1+\alpha}}$ in formula~(\ref{eval1})
after the substitution~(\ref{lemma_2});
then it is easy to see that the ratio
$$
\frac{\cN(\cM \setminus\cA)}{\cN(\cM)}\approx \frac 1{N^m},
$$
where $m$ is an arbitrary integer,
holds for $\Delta=N^{3/4+\varepsilon}$.
The proof of the theorem is complete.

We prove a cumulative formula in which the densities coincide in
shape with the Bose--Einstein distribution. The difference
consists only in that, instead of the set $ x_n$ of random
variables or eigenvalues of the Hamiltonian operator, the formula
contains some of their averages over the cells. In view of our
theorem, the $\varepsilon_i$, which are averages of the energy $
x_k$ at the $i$th cell, are nonlinear averages in the sense of
Kolmogorov~\cite{NelinSred}.

In conclusion, I express my deep thanks
to G.~V.~Koval' and V.~E.~Nazaikinskii
for fruitful discussions.

\end{document}